\documentclass{article}
\usepackage{amssymb}
\usepackage{amsmath}
\input epsf.sty

\begin{document}

\def\phi{{\varphi}}
\def\d{{\rm d}}
\def\D{{\rm D}}
\def\Q{{\mathbb Q}}
\def\C{{\mathbb C}}
\def\N{{\mathbb N}}
\def\P{{\mathbb P}}
\def\R{{\mathbb R}}
\def\deg{{\rm deg\,}}
\def\Re{{\rm Re\,}}
\def\Det{{\rm Det}}
\def\dim{{\rm dim\,}}
\def\Gal{{\rm Gal\,}}
\def\PGL{{\rm PGL\,}}
\def\Rat{{\rm Rat\,}}
\def\Aut{{\rm Aut\,}}
\def\St{{\rm St\,}}

\title
{Conservative polynomials and yet another
action of $\Gal(\bar \Q/\Q)$
on plane trees}
\author{F. Pakovich}

\date{}
\maketitle

\def\bp{\begin{proposition}}
\def\ep{\end{proposition}}
\def\bt{\begin{theorem}}
\def\et{\end{theorem}}
\def\be{\begin{equation}}
\def\l{\label}
\def\ee{\end{equation}}
\def\bl{\begin{lemma}}
\def\el{\end{lemma}}
\def\bpr{\begin{problem}}
\def\epr{\end{problem}}
\def\bc{\begin{corollary}}
\def\ec{\end{corollary}}
\def\pr{\noindent{\it Proof. }}
\def\note{\noindent{\bf Note. }}
\def\bd{\begin{definition}}
\def\ed{\end{definition}}
\newtheorem{theorem}{Theorem}[section]
\newtheorem{lemma}{Lemma}[section]
\newtheorem{definition}{Definition}[section]
\newtheorem{corollary}{Corollary}[section]
\newtheorem{proposition}{Proposition}[section]
\newtheorem{problem}{Problem}[section]

\section{Introduction} The remarkable Belyi theorem \cite{be} states that an algebraic curve $X$ defined over $\C$ 
is defined over $\bar\Q$ if and only if there is a holomorhic function $\beta\, :\, X \rightarrow \C\P^1$ which is ramified only over $\{0,1,\infty\}.$ A pair consisting
of such $X$ and $\beta$ is called a Belyi pair. The isomorphism classes of Belyi pairs
are in one-to-one correspondence with the classes of isotopical equivalence of bicolored graphs embedded into topological models of Riemann surfaces. Namely, if $(X,\beta)$
is a Belyi pair then the corresponding graph $\Omega$ is the preimage of the segment $[0,1]$ under the map $\beta\, :\, X \rightarrow \C\P^1$ where white (resp. black) vertices 
of $\Omega$ are preimages of $0$ (resp. $1$). 

The absolute Galois group $\Gamma=\Gal(\bar \Q/\Q)$ acts naturally on isomorphism classes of Belyi pairs and this action descends to an action 
on equivalence classes of bicolored graphs. The study of this action, which we will denote by $G$, is the subject of the Grotendieck theory
of ``Dessins d'enfants'' (see e.g. \cite{lz} and the bibliography there). Note that the action $G$ is highly non-trivial even when restricted on the bicolored plane trees 
(which are in one-to-one correspondence with the equivalence classes of
the polynomial Belyi functions on the Riemann sphere). In particular, the
action of $\Gamma$ on plane bicolored trees is faithful (\cite{sh}).

The Belyi pairs are an example of ``rigid'' analytical objects which are completely defined up to an equivalence by some ``combinatorial'' data. Another example of such objects  
is postcritically finite rational functions  
that is the rational functions for which orbits of critical points under iterations are finite. Indeed, by the result due to Thurston \cite{dh} such functions, apart from a very special family, 
are uniquely defined,
up to a conjugacy by a M\"obius transformation, by purely combinatorial data. In particular, conjugacy classes of postcritically finite polynomials can be classified by means some rather complicated
combinatorial objects called Hubbard trees (see \cite{poi1}, \cite{poi2}, \cite{pil}).

The finiteness of conjugacy classes of postcritically finite polynomials of a given degree 
makes possible to define an action of $\Gamma$ on these classes similarly to the action $G.$
Note that this action, which we will denote by $D$, is faithful \cite{pil}. Furthermore, there is some interplay between two theories since each equivalence class of polynomial
Belyi functions  
contains postcritically finite polynomials
(see \cite{pil}).

The simplest example of postcritically finite polynomials is polynomials
with all critical points fixed, called conservative,
and in this paper we study 
the action $D$ restricted on such polynomials.
This particular case seems to be especially interesting since,
in distinction with the general case,
for conservative polynomials the corresponding combinatorial data can be described quite transparently. Namely, by the result due to Tischler \cite{t} 
the
equivalence classes of conservative polynomials as holomorhic dynamical systems on $\C$ are in one-to-one correspondence with the classes of isotopical equivalence of bicolored plane trees. So, the 
absolute Galois group $\Gamma$ acts on bicolored plane trees in two different ways:
one action is induced by the action on the polynomial Belyi functions and the other one
by the action on the conservative polynomials !

The paper is organized as follows. In the second section we reproduce the Tishler 
correspondence between conservative polynomials and bicolored plane trees
and provide some examples.
We write explicitly equations for defining a conservative polynomial $C(z)$ corresponding to a tree $\lambda$ and show how $C(z)$ reflects the symmetries of $\lambda.$
Besides, we prove that 
any conservative polynomial is 
indecomposable (that is can not be represented as a composition of two non-linear polynomials) unless it is equivalent to the polynomial $z^n$ for composite $n\in \N.$

In the third section we establish some properties of the action $D$ 
and compare it with the action $G$. 
In particular, 
we show that the list of valencies of ``white" (but not ``black" !) vertices of a tree $\lambda$ 
and the symmetry group of $\lambda$
are combinatorial Galois invariants of the action $D$.
Then we give several examples of calculations 
of Galois orbits.
Finally, we describe all trees which are defined uniquely by the list 
of valencies of ``white" vertices 
and calculate the 
corresponding polynomials. 




\section{Conservative polynomials and plane trees
} 

Recall that {\it a plane tree} is a tree embedded into the plane and that {\it a bicolored tree} is a tree vertices of which are colored in two colors
in such a way that 
any edge connects vertices 
of different colors.
Two bicolored plane trees $\lambda$ and $\tilde \lambda$ are called {\it equivalent} if 
there exists an orientation preserving homeomorphism $h$ of $\C$ such that 
$\tilde \lambda=h(\lambda)$ and $h$ preserves the 
colors of vertices.

A complex polynomial $C(z)$ is called {\it conservative} if all
its critical points are fixed that is if the equality $C^{\prime}(\zeta)=0,$ $\zeta\in \C,$ 
implies that $C(\zeta)=\zeta.$ A conservative polynomials $C(z)$ is called {\it normalised}
if $C(z)$ is monic and $C(0)=0.$ 
Two conservative polynomials $\tilde C(z)$, $C(z)$ are called {\it equivalent} if 
there exists a complex polynomial $A(z)$ of degree one
such that $\tilde C=A^{-1}\circ C\circ A.$ 

Conservative polynomials were introduced by
Smale \cite{s} in connection with his ``mean value conjecture".
Motivated by Smale's conjecture Kostrikin proposed in \cite{k} 
several conjectures concerning conservative polynomials.
In particular, on the base of 
numerical experiments Kostrikin conjectured that the
number of normalised conservative polynomials of degree $d$
is finite and is equal to $C_{2d-2}^{d-1}.$
This
conjecture was proved by Tischler in paper \cite{t}. Moreover, in
this paper a one-to-one correspondence between equivalence classes of
conservative polynomials and
equivalence classes of plane bicolored trees was
established. Below we reproduce the Tishler construction.

Let $\zeta$ be a critical point of $C(z)$ and $d\geq 2$ be the multiplicity
of $C(z)$ at $\zeta.$ Then one can show (see \cite{t}) that the immediate attractive
basin $B_{\zeta}$ of $\zeta$ is a disk and that
there is an analytic conjugation of $C(z)$ on $B_{\zeta}$ 
to $z\rightarrow z^{d}$ on the unit disk $D$ such that the conjugating map $\phi_{\zeta} : \D \rightarrow
B_{\zeta}$ extends continuously to the closed unit
disk $\bar \D$. Let $S$ be a union of $d-1$ radial
segments which are forward invariant under the map $z\rightarrow z^{d}$
on $\bar \D$ and $S_{\zeta}$ be the image of $S$ under the map $\phi$
(see Fig. 1 where $d=4$). We consider $S_{\zeta}$ as a bicolored graph 
with a unique white vertex which is the image of zero and $d-1$ 
black vertices which are the images of end-points of $S$.

\vskip 0.2cm
\medskip
\epsfxsize=8.5truecm
\centerline{\epsffile{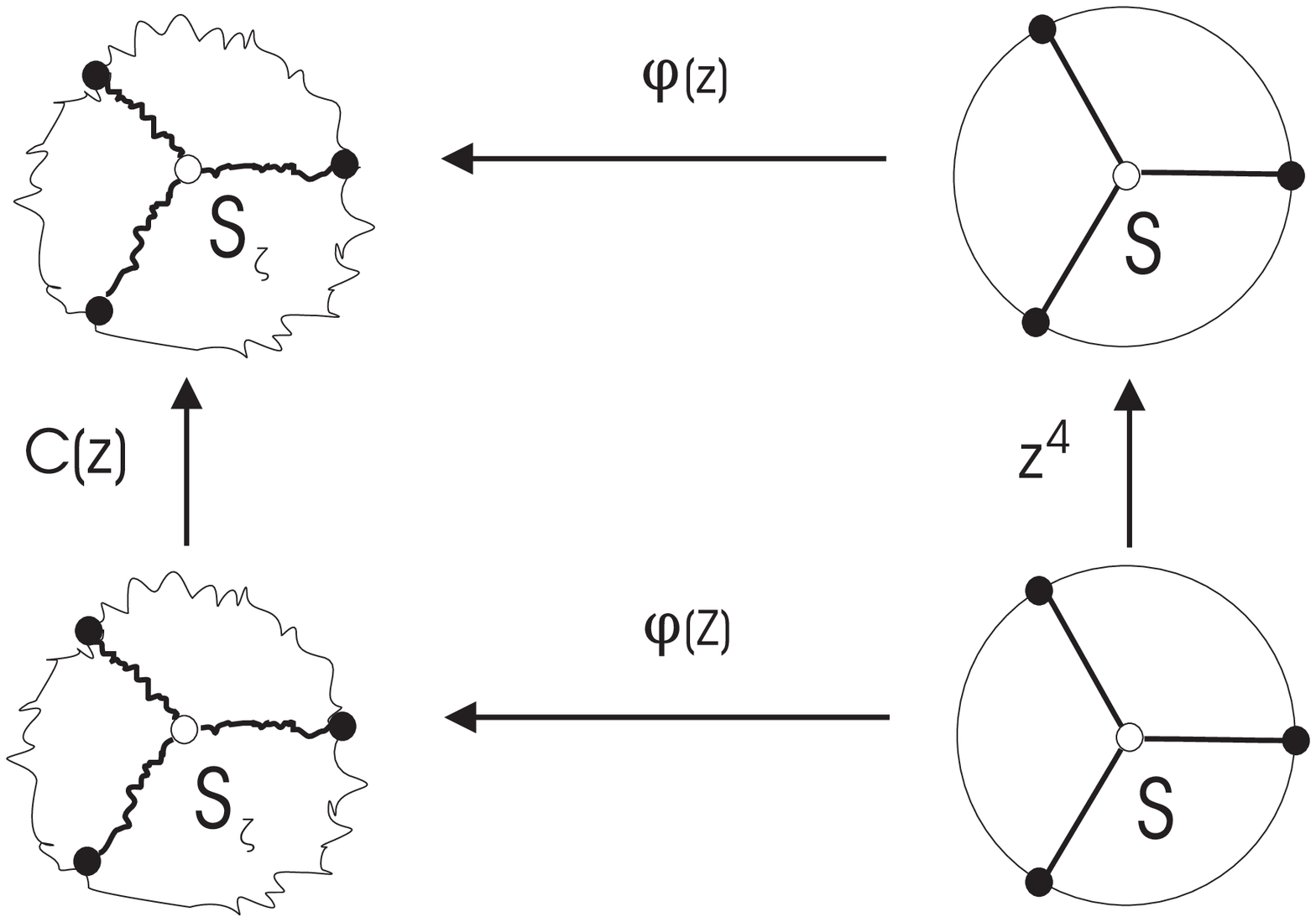}}
\smallskip
\centerline{Figure 1.}
\vskip 0.2cm

Define now a bicolored graph $\lambda_C$ as a union
$\lambda_C=\cup_{i=1}^p S_{\zeta_i},$ where $\zeta_i,$ \linebreak $1\leq i \leq p,$ are all 
finite critical points of $C(z).$

Clearly,
the valency of a white vertex $v$ of $\lambda_C$ coincides with the multiplicity of the point $v$ with respect to the map $C^{\prime}(z): \mathbb C \rightarrow \mathbb C.$ 
Note also that by construction the graph $\lambda_C$ is a forward 
invariant of $C(z)$ and white (resp. black) vertices of $\lambda_C$ are attractive (resp. repelling) fixed points of $C(z).$  

It turns out that the graph
$\lambda_C$ is actually a tree. Moreover, the following theorem
proved by Tischler is true. 

\bt [\cite{t}]
The map $C\rightarrow \lambda_C $
descends
to a bijection between equivalence classes of conservative
polynomials of degree $d$ and
equivalence classes of bicolored plane trees with $d-1$ edges. Furthermore, the number
of normalised conservative polynomials of degree $d$ is $\binom{2d-2}{d-1}.$

\et

The simplest example of a conservative polynomial is the polynomial
$z^d$ with a unique critical point $\zeta=0.$
Clearly, the corresponding tree is a $d-1$-edged
star.

Another 
example (cf. \cite{s},\cite{k}) is the polynomial $f_d(z)=z^d+(d/d-1)z,$ $d\geq
2$. 
Since all 
zeros 
of $f_d^{\prime}(z)$ are simple all white vertices of the 
corresponding tree $\lambda_d$ are of valency 1 and
therefore $\lambda_d$ is also a
$d-1$-edged star but {\it with the
bicoloring changed.} 
On the right side of Fig. 2 the dynamics of $f_6(z)$ are shown: the points for which 
the limits of iterations have the same value 
are painted by the same color\footnote{To prepare dynamical pictures for this paper we used C. McMullen's programs available on http://www.math.harvard.edu/~ctm/programs.html}. 

Note that we face here a phenomenon which is absent in the
``Dessins d'enfants" theory: the calculation of a polynomial corresponding to a tree
$\lambda$ and the calculation of a polynomial corresponding to a tree $\lambda^*$ obtained
from $\lambda$ by the change of the bicoloring are essentially different questions apparently not connected 
between themselves.

\vskip 0.2cm
\medskip
\epsfxsize=9truecm
\centerline{\epsffile{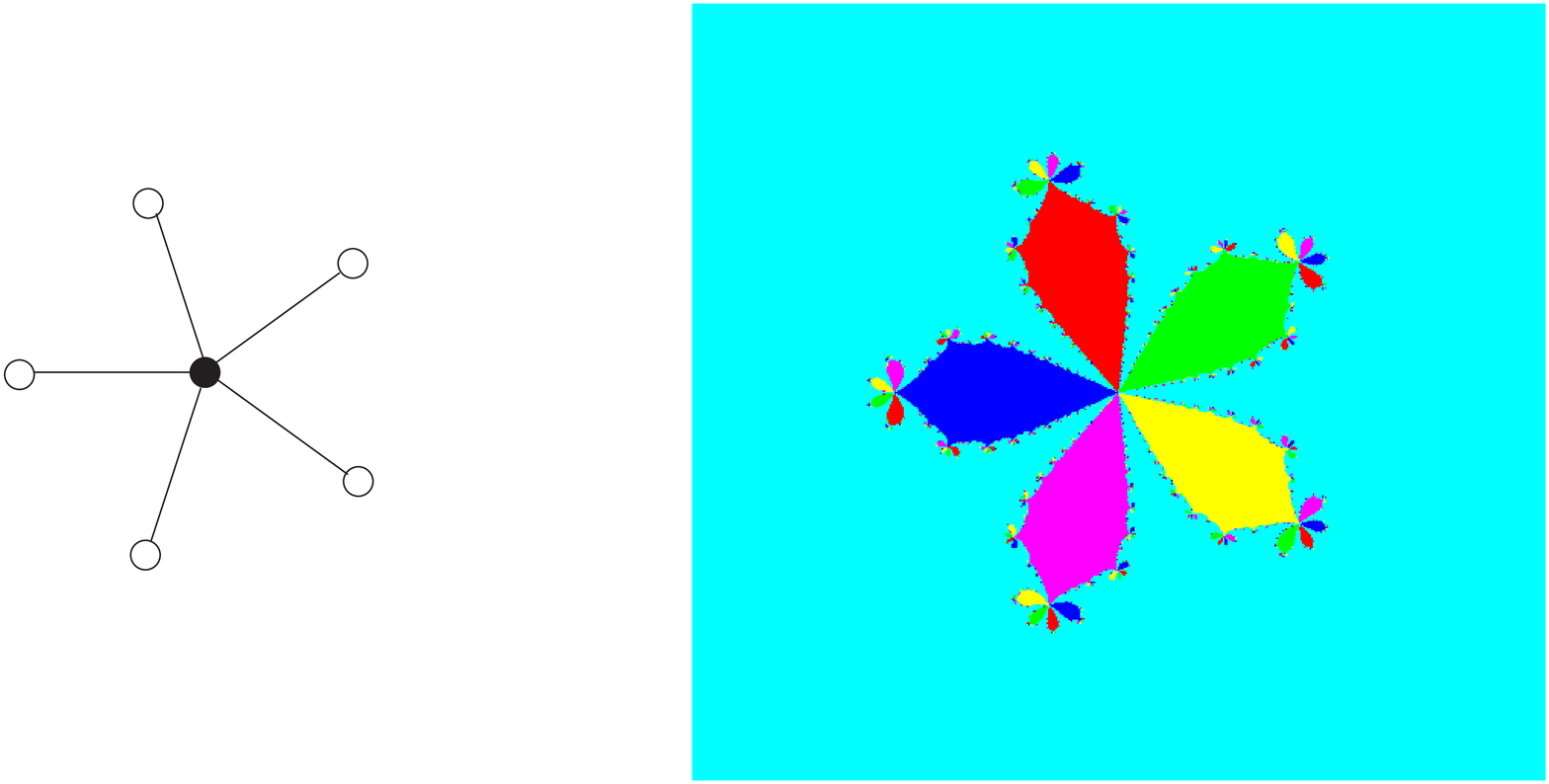}}
\smallskip
\centerline{Figure 2.}
\vskip 0.2cm

Starting from a tree $\lambda$ we can find a conservative polynomial 
$C(z)$ from the corresponding equivalence class as follows.  
Let $\alpha=<\alpha_1,\alpha_2, ... \,,\alpha_p>$ be 
the sequence of valencies of white vertices of $\lambda$ in decreasing order.
Denote by $a_1, a_2, ... , a_p$ unknown
coordinates of white vertices of $\lambda$ 
coinciding with zeros of $C^{\prime}(z).$ 
Clearly, without loss of generality we can assume that $a_1=0.$

Set $C(z)$ equal to the indefinite integral $$\int nz^{\alpha_1}(z-a_2)^{\alpha_2}\, ...\,\, (z-a_{p})^{\alpha_p}\d z$$ normalised by the condition $C(0)=0.$  
Then the system of equations to determine $a_2, a_3, ... , a_{p}$ is 
\be \l{ss} C(a_j)=a_j,\ \ \ \ \ \ \ \ \ \ 2\leq j \leq p. \ee

Observe that system \eqref{ss} depends only on the sequence $\alpha$ which is called {\it the type} of $\lambda$. Therefore, solutions of \eqref{ss} along with a polynomial corresponding to $\lambda$ contain any polynomial 
corresponding to a tree of type $\alpha$. Note also that the system \eqref{ss}
may have solutions for
which the numbers $a_1, a_2, ...\, ,a_{p}$ are {\it not} mutually distinct.
Such a solution also corresponds to a tree $\tilde \lambda$ but the type of this tree is distinct from $\alpha.$ Geometrically, $\tilde \lambda$ is obtained from some tree $\lambda$ 
of type $\alpha$ by ``merging'' some number of white vertices of $\lambda$.

Let us describe now following \cite{t} the combinatorial datum which permits to determine a normalised conservative 
polynomial corresponding to a tree $\lambda$ uniquely. First, the condition $f(0)=0$ corresponds to the choice of a vertex $v$ of $\lambda$ which we place at the origin. Furthermore, since $C(az)/a,$ $a\in \C$ is monic if and only if $a^{d-1}=1,$ 
where $d=\deg C(z),$ after choosing $v$ the corresponding normalised polynomial 
is defined up to a change $z\rightarrow \varepsilon z$, $\varepsilon^{d-1}=1.$  

Let $S_{\infty}\subset \C\P^1$ be a graph defined like the graphs
$S_{\zeta_i},$ $1\leq i \leq p,$ with only difference that we start from the infinite critical point of $C(z).$ 
It is shown in \cite{t} that each edge of $S_{\infty}$ ends at some black vertex of $\lambda_C$ and that the number of edges ending at a black vertex $w$ is equal to the multiplicity of $w$. More precisely, each angle formed by two adjacent edges emerging from 
a black vertex of $\lambda$ contains exactly one edge of $S_{\infty}$ (see Fig. 3). 

\vskip 0.2cm
\medskip
\epsfxsize=6truecm
\centerline{\epsffile{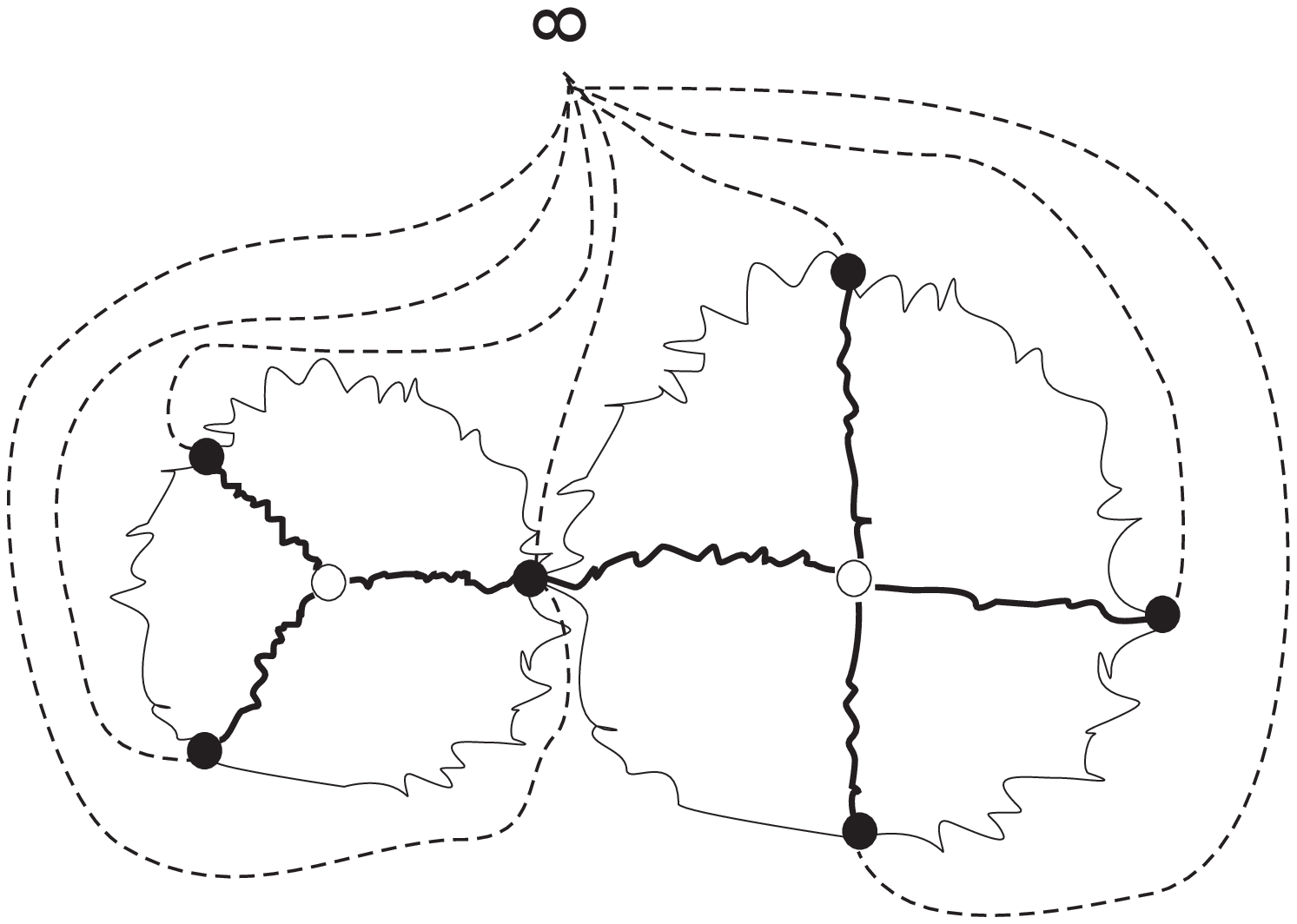}}
\smallskip
\centerline{Figure 3.}
\vskip 0.2cm

Clearly, a choice of a $d-1$-th root of unity
corresponds to a choice of an edge of $S_{\infty}$ and therefore 
a normalised conservative 
polynomial corresponding to a tree $\lambda$ is determined uniquely by fixing
a vertex $v$ (black or white) of $\lambda$ and an ``angle'' $\phi$
adjacent to a black vertex of $\lambda$.

As an example consider the set of bicolored trees with $4$ vertices. It is easy to see that there are exactly three
such trees which are shown on Fig. 4.

\vskip 0.2cm
\medskip
\epsfxsize=9truecm
\centerline{\epsffile{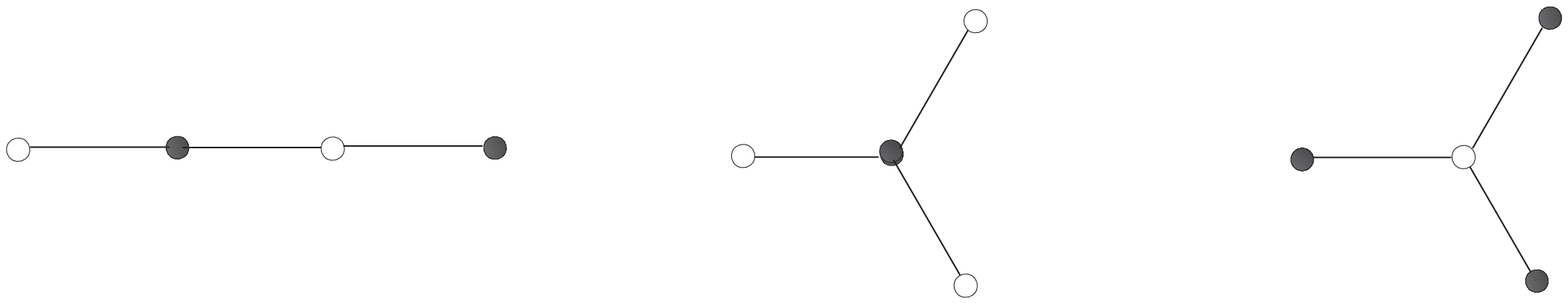}}
\smallskip
\centerline{Figure 4.}
\vskip 0.2cm
\noindent On the other hand, by the Tishler theorem there exist $\binom{6}{3}=20$ different normalised conservative polynomials of degree $4:$ to the first tree shown on Fig. 4 correspond 12 different polynomials
while to the second and to the third ones correspond only 4 different polynomials.

In general, since collections $(\lambda, v_1, \phi_1)$ and $(\lambda, v_2, \phi_2)$ are equivalent 
if and only if there exists an orientation preserving homeomorphism $\chi$ of $\C$ 
such that $\chi(\lambda)=\lambda$ and $\chi(v_1)=v_2,$  $\chi(\phi_1)=\phi_2,$ it is easy 
to see that 
to a plane tree $\lambda$ with $d$ vertices correspond $d(d-1)/\vert \Aut \lambda \vert$ different normalised conservative polynomial, where $\Aut \lambda$ denotes
the group of symmetries of $\lambda.$ Furthermore, the following statement holds.

\bp \l{x} A tree $\lambda$ has a symmetry of order $k$ if and only if the corresponding class
of conservative polynomials contains a polynomial of the form 
$C(z)=zR(z^k),$ where $R(z)$ is a polynomial. 
\ep

\pr Indeed, 
if  
$C(z)$ is a normalised conservative polynomial of degree $d$ corresponding to a collection $(\lambda, v, \phi)$ and $v$ is placed at the origin then to a collection  
$(\lambda, v, \tilde \phi)$ corresponds a conservative polynomial $C(\varepsilon z)/\varepsilon,$ where $\varepsilon$ 
is some $d-1$-th root of unity. Therefore, a tree $\lambda$ has a
symmetry of order $k$ with the center at $v$ if and only if $C(z)=C(\varepsilon z)/\varepsilon$ for 
any $k$-th root of unity $\varepsilon$. This condition is equivalent to the condition that the polynomial $C(z)/z$ is invariant with respect
to any rotation of the form $z\rightarrow \varepsilon z,$ where $\varepsilon$ is a 
$k$-th root of unity. In its turn the last condition is equivalent to the condition that $C(z)=zR(z^k)$ for some polynomial $R(z).$

\vskip 0.2cm 
Since many constructions of the ``Dessins d'enfants" theory
(for instance, the Belyi theorem) make use compositions of functions
and such compositions are survive under the Galois action
it is natural to ask about compositional properties of the conservative polynomials.
It turns out that conservative polynomials are essentially indecomposable.

\bp\l{prim}{All conservative polynomials not equivalent to $z^n$ 
for composite $n$ are indecomposable.} 
\ep
\pr
Indeed, suppose that 
$C(z)=C_1(C_2(z))$ with $\deg C_1(z),\deg C_2(z) >1.$ Let $\zeta\in \C$ be a   
critical point of the polynomial $C_1(z).$ Then the  chain rule implies 
that any point $\mu\in \C$ such that $C_2(\mu)=\zeta$ is a critical 
point of $C(z).$ Furthermore, if the polynomial $C_2(z)$ is not equal to $A(z-\mu)^l+\zeta$ 
for some $A,\mu, \zeta \in \C$ and integer $l\geq 2,$ then there exist $\mu_1,\mu_2\in \C,$ $\mu_1\neq \mu_2,$ such that $C_2(\mu_1)=C_2(\mu_2)=\zeta.$ For these points we have 
$C^{\prime}(\mu_1)= C^{\prime}(\mu_2)=0$ and $C(\mu_1)= C(\mu_2)
=C_1(\zeta).$ Since 
$\mu_1\neq \mu_2$ this contradicts to the condition that $\mu_1,\mu_2$ are 
fixed points of $C(z).$ 

Therefore, $C_2(z)=A(z-\mu)^l+\zeta.$ In particular, for any $\tilde\zeta\in \C,$ $\tilde\zeta\neq \zeta$   
there exist 
mutually distinct $\tilde\mu_1,\tilde\mu_2, \,...\, , \tilde\mu_l\in \C$ such that 
$C_2(\tilde\mu_i)=\tilde\zeta,$ $1\leq i \leq l.$ Thus if $\tilde\zeta$ is 
a critical point of $C_1(z)$ distinct from $\zeta$ then as above we obtain a contradiction since
$C^{\prime}(\tilde \mu_1)= C^{\prime}(\tilde \mu_2)=0$ and $C(\tilde \mu_1)= C(\tilde \mu_2)=C_1(\tilde \zeta).$ Hence, $\zeta$ is a unique finite critical point of
$C_1(z).$ It follows that $C_1(z)=B(z-\zeta)^{k}+\nu$ for some $B,\nu \in \C$ and integer $k\geq 2$ and hence
$C(z)=C_1(C_2(z))$ is equivalent to $z^n$ 
for composite $n.$    
\vskip 0.2cm


\section{Galois group action}
In each equivalence class of conservative polynomials there exist 
polynomials with algebraic coefficients. Indeed, it follows from the
Tischler
theorem that system \eqref{ss} has only a finite number of solutions.
Therefore, since equations \eqref{ss}
have rational coefficients, 
all these solutions are algebraic. 

Furthermore, the group $\Gamma$ acts on the set of conservative 
polynomials with algebraic coefficients in a natural way: it is easy to see that
if $C(z)$ is a conservative polynomial with algebraic
coefficients and $\sigma \in \Gamma$   
then the polynomial $C^{\sigma}(z)$ obtained from $C(z)$ by the action of $\sigma$ on 
coefficients of $C(z)$  
again is
a conservative polynomial.  
Moreover, since $(A^{-1}\circ P\circ A)^{\sigma}=(A^{\sigma})^{-1}\circ P^{\sigma}\circ A^{\sigma},$ this action descends to an action on equivalence classes. 
Hence, by the Tischler theorem, we obtain an action $D$ of $\Gamma$ on bicolored plane trees. 

The type $\alpha$ of a tree $\lambda$ is an invariant of the action $D$ since the sequence $\alpha$ coincides with the sequence of multiplicities of zeros 
of $C^{\prime}(z).$ Furthermore, like the ``Dessins d'enfants" theory the Galois
orbit of $\lambda$ often coincides with the set of all trees of 
type $\alpha$. Moreover, in view of proposition \ref{x} the symmetry group of $\lambda$
is also a Galois invariant of the action $D$ since if 
$C(z)$ has the form $C(z)=zR(z^k)$ for some polynomial $R(z)$ then the Galois conjugated polynomial also has such a form. 

Finally, it is easy to see (see e. g. Example 1 
below) that the action $D$ is distinct from the action $G.$ 
Summing up we obtain the following statement.

\bt \l{y} $D$ is a well-defined action of the group $\Gal(\bar \Q/\Q)$ on
bicolored plane trees distinct from the action $G.$ The type and the symmetry group of a tree $\lambda$ are combinatorial Galois invariants of the action $D.$
\et

Define {\it a field of modules} $k^D_{\lambda}$ of a tree $\lambda$ as a fixed field of the stabilizer of
$\lambda$ with respect to the action $D.$ Since there exist only a finite number 
of trees of a given type the field $k^D_{\lambda}$ is a number field
whose degree over $\Q$ is equal to the
length of the orbit containing $\lambda.$ Note that the result proved in \cite{sil} implies
that for any bicolored
plane tree $\lambda$ there exists a conservative polynomial $C(z)$ 
from the corresponding equivalence class 
such that $C(z)\in k_{\lambda}^D[z].$ 

\vskip 0.2cm
\noindent{\bf Example 1.} As the first example consider the trees of type $<3,1,1>$.
There exist two trees of this type which are shown
on Fig. 5.

\medskip
\epsfxsize=7truecm
\centerline{\epsffile{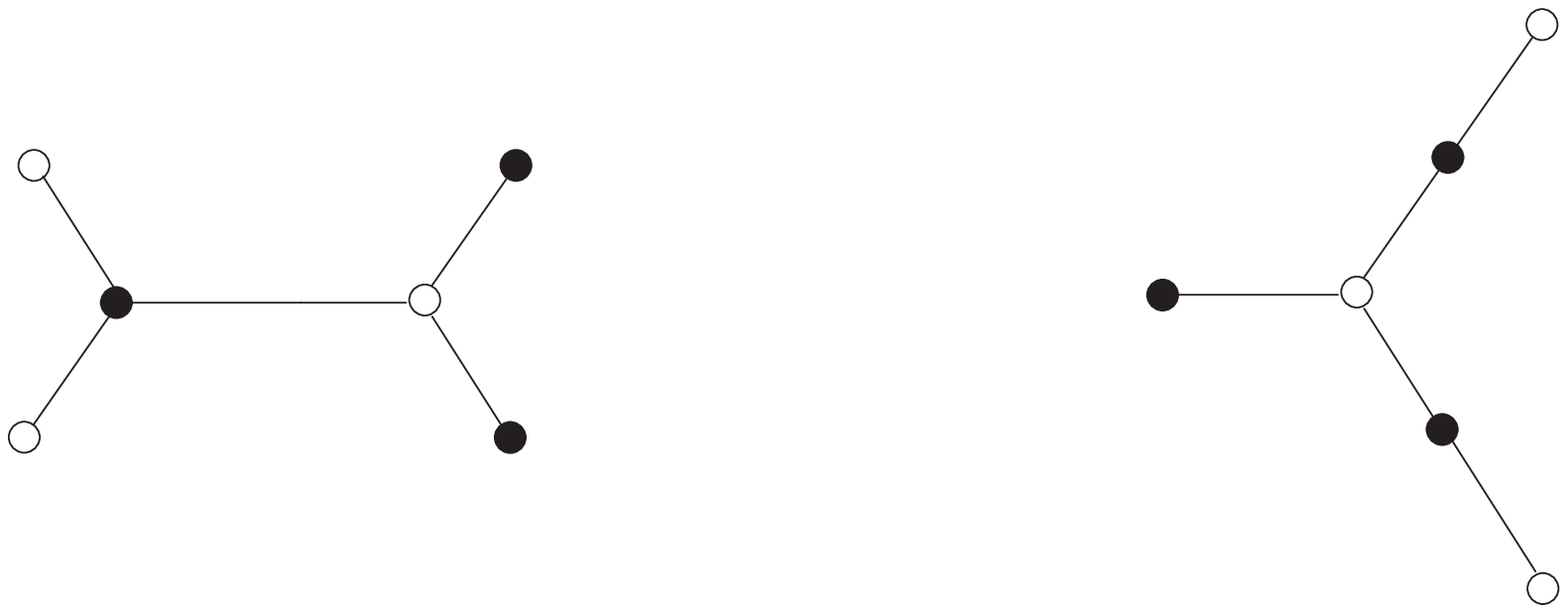}}
\smallskip
\centerline{Figure 5.}
\vskip 0.2cm 
\noindent Place the white vertex of valency 3 at zero. Then 
$$C^{\prime}(z)=az^3(z^2+cz+b) $$
for some $a,b,c \in \C$ such that the polynomial $z^2+cz+b$ has two different roots distinct from 0. Furthermore, 
$c\neq 0$ since otherwise 
$C(z)=\tilde C(z^2)$ for some polynomial $\tilde C(z)$ in contradiction with
proposition \ref{prim}. Therefore, we can suppose that $c=1.$ 

Since $C(0)=0$ we have  
$$C(z)=a(z^6/6+z^5/5+bz^4/4)$$ and a calculation shows that
\be\l{rys}
C(z)-z=A(z)(z^2+z+b)+B(z),
\ee
where 
\be
\begin{split}
A(z)&=az^4/6+az^3/30+(ab/12-a/30)z^2+(-7ab/60+a/30)z \\
&-ab^2/12+3ab/20-a/30,
\\
B(z)&=(a/30-1+ab^2/5-11ab/60)z+ab^3/12-3ab^2/20+ab/30.
\notag
\end{split}
\ee 

Since roots $\beta_1, \beta_2$ of $z^2+z+b$ are fixed points of $C(z)$ and $\beta_1\neq \beta_2$ it follows from 
\eqref{rys} taking into account the equality $\deg B(z)=1$ 
that $B(z)\equiv 0.$ This gives us the system    

$$
\left\{
\begin{aligned}
a/30-1+ab^2/5-11ab/60&=0,\\
ab^3/12-3ab^2/20+ab/30&=0.\\
\end{aligned}
\right. 
$$ 
Solving this system we conclude that either
\be \l{sol} a=30, \ \ \ \ b=0,\ee or $$a=-5055/8 \pm 795\sqrt{41}/8,
\ \ \ \ b=9/10 \pm \sqrt{41}/10.$$ 
All three solutions above correspond some conservative polynomials but in case
when \eqref{sol} holds  
one of the roots of $z^2+dz+b$ coincides with zero.
This means that solution \eqref{sol} actually corresponds to the tree $\lambda_{4,1}$ (see Fig. 8 below).

After rejecting solution \eqref{sol} it remains two Galois conjugated solutions 
and therefore the trees shown on Fig. 5 form a two-element Galois
orbit with field of modulus $\Q(\sqrt{41}).$
\vskip 0.2cm  

In this example we can see an important difference between the actions $G$ and $D.$ 
Namely, 
for the action $G$ not only the 
list of valencies of ``white" vertices but also the 
list of valencies of ``black" vertices is Galois invariant. Therefore, 
``expected'' orbits of the action $D$ 
are much longer than the ones with respect to the action $G$.
In particular, since both of trees shown on Fig. 5
are determined uniquely by the lists of valencies of ``white" and ``black" vertices,
each of them forms an one-element Galois orbit
with respect to the action $G.$

\vskip 0.2cm

\noindent{\bf Example 2.} Consider now the trees of type $<3,1,1>$.
Again there exist two trees of this type: 
one of them is shown on Fig. 6 and the other one on Fig. 7. 
Nevertheless, since the first tree has a symmetry of order 2 and the second one does not have symmetries
it follows from theorem \ref{y} that each of these trees form an one-element Galois orbits
with respect to the action $G.$

In order to calculate a conservative polynomial corresponding to the first tree set 
$$C^{\prime}(z)=az^2(z^2+bz-1).$$ In other words we
suppose that the white vertex of valency 2 is located at zero and that
the product of the coordinates of two white vertices of valency 1
equals -1. Observe now that in view of proposition \ref{x} we necessarily have $b=0.$
In particular, the coordinates of white vertices of valency 1
are $\pm 1.$
 
Therefore, $C(z)=az^5/5-az^3/3$ and the
conditions $C(1)=1,$ $C(-1)=-1$ reduce to the equality $-2a/15=1.$
Hence,  
$$C(z)=-3z^5/2+5z^3/2.$$ On the right side 
of Fig. 6 the
dynamics of $C(z)$ are shown.

\vskip 0.2cm
\medskip
\epsfxsize=10truecm
\centerline{\epsffile{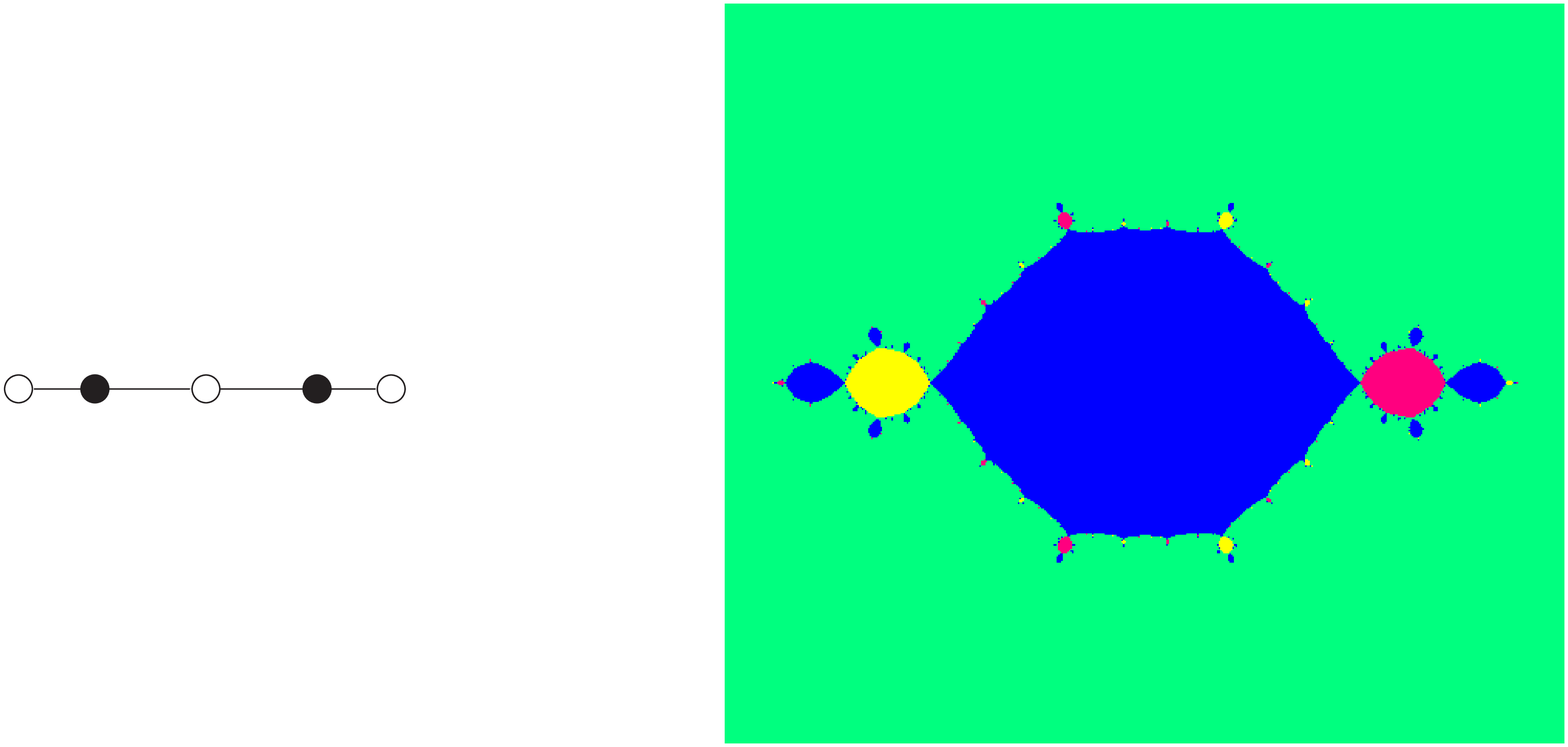}}
\smallskip
\centerline{Figure 6.}
\vskip 0.2cm
\vskip 0.2cm
To calculate a polynomial corresponding to the second tree observe that in this case 
proposition \ref{x} implies that the sum of coordinates of two white vertices of valency 1 is necessarily
distinct from zero so we can set this sum  
equal 2. 


Then 
$$C^{\prime}(z)=az^2(z^2+2z+b), \ \ \ \ C(z)=a(z^5/5+z^4/2+bz^3/3)$$ 
and as in the first example the remainder 
$$B(z)=(14ab/15-4a/5-1-2ab^2/15)z+11ab^2/30-2ab/5 $$ after division of
of $C(z)-z$ by $z^2+2z+b$ equals zero. Solving the
corresponding system 
\be \l{ls}
\left\{
\begin{aligned}
14ab/15-4a/5-1-2ab^2/15&=0,\\
11ab^2/30-2ab/5&=0\\
\end{aligned}
\right. 
\ee
and rejecting the solution $$a=-5/4, \ \ \ \ b=0$$ 
corresponding to the tree $\lambda_{2,1}$ we conclude that $$a=605/36, \ \ \ \ b=12/11$$ and therefore 
$$C(z)= 121z^5/36+605z^4/72+55z^3/9.$$ 
\vskip 0.2cm
\medskip
\epsfxsize=10truecm
\centerline{\epsffile{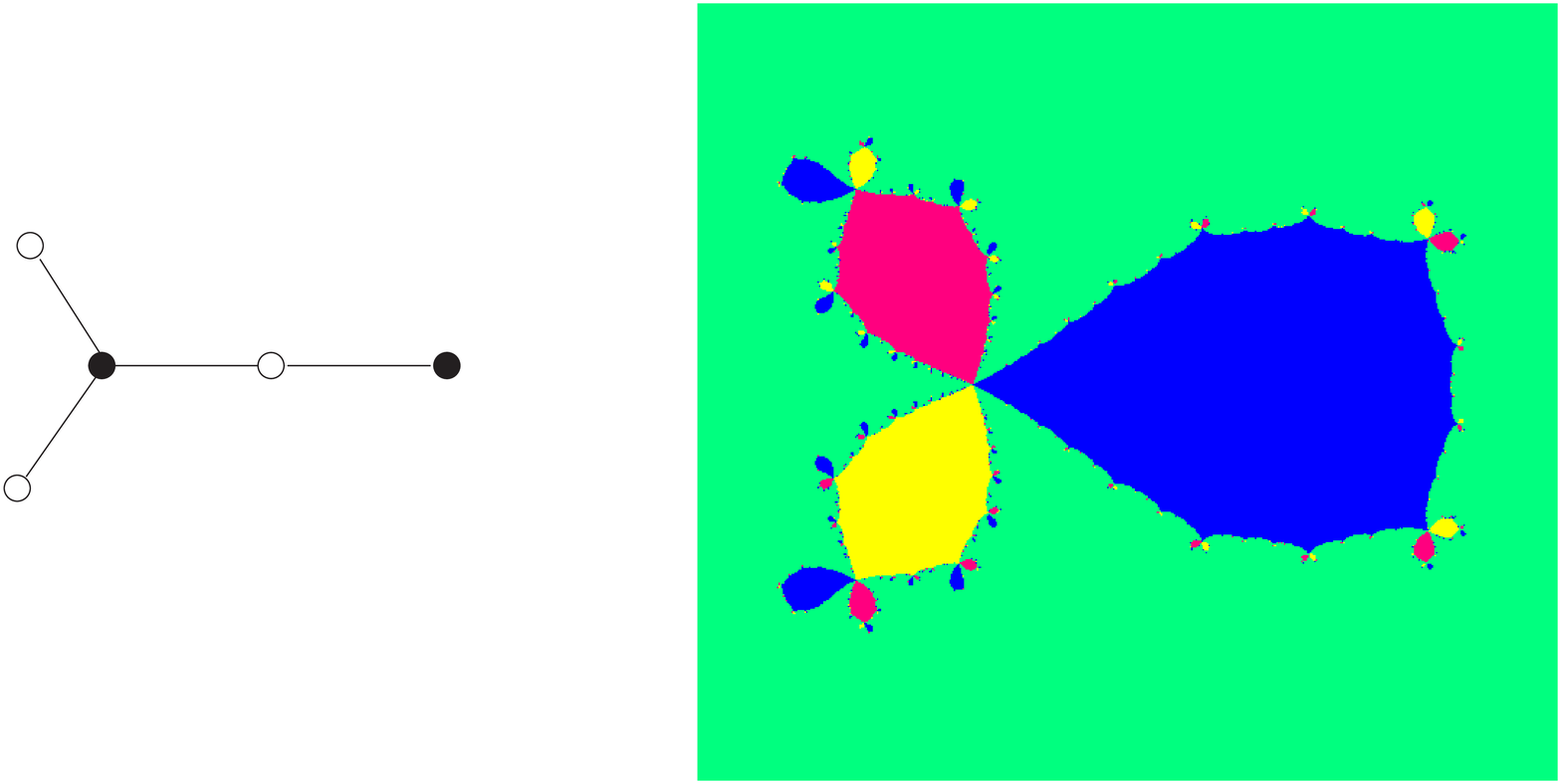}}
\smallskip
\centerline{Figure 7.}
\vskip 0.2cm
\vskip 0.2cm

\noindent{\bf Example 3.} Consider finally the tree $\lambda_{r,s}$
shown on Fig. 8, where $r,s$ are some integers $\geq 1$.

\vskip 0.2cm
\medskip
\epsfxsize=4truecm
\centerline{\epsffile{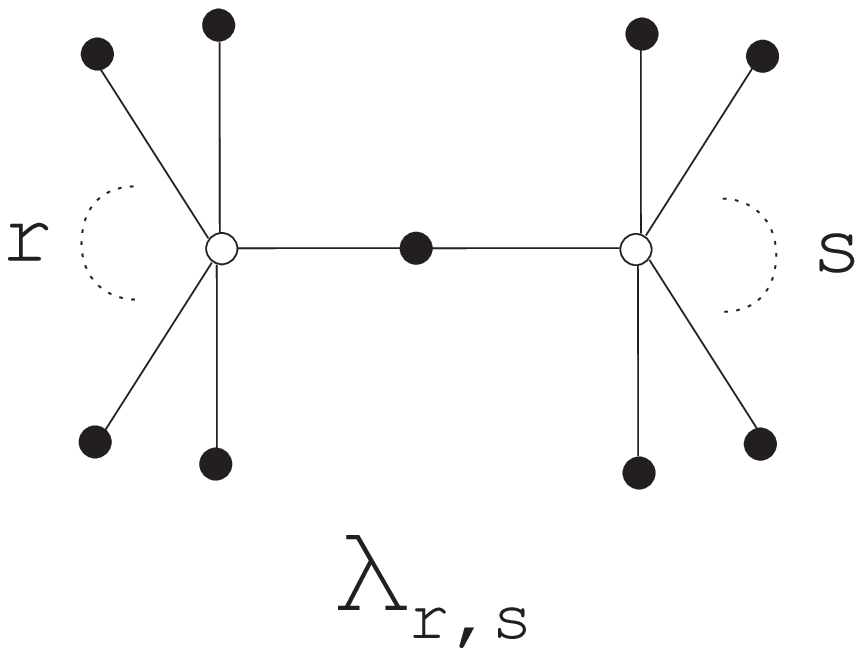}}
\smallskip
\centerline{Figure 8.}
\vskip 0.2cm

\noindent Clearly, we can place white vertices to the points $0,1.$ 
Then
$$C'(z)=cz^r(1-z)^s,$$
where $c\in \C$ is a parameter to define. Integrating and taking into account that $C(0)=0,$ we conclude that
$$C(z)=cz^{r+1}\ _2F_1(-s,r+1,2+r,z). $$ Finally, 
since 
$$C(1)=
c\frac{r!s!}{(r+s+1)!},$$
the condition $C(1)=1$ implies that $$c=\frac{(r+s+1)!}{r!s!}.$$ 

\vskip 0.2cm 

Note that any tree $\lambda_{r,s}$ is determined uniquely by its type. Since in view of 
theorem \ref{y} such a tree necessarily forms an one-element Galois orbit it is 
interesting to know how many trees possess
this property.
The proposition below
shows that such a phenomenon is rather exceptional.

\bp Suppose that a tree $\lambda$ is determined uniquely by its type. Then $\lambda$ is either a star with
black center, or a star with
white center, or the tree $\lambda_{r,s}$ for some $r,s \geq 1.$
\ep

\pr 
Let $\alpha_1,\alpha_2, ... ,\alpha_p$ (resp. $\beta_1,\beta_2, ... ,
\beta_q$) be the list of valencies of white (resp. black) vertices
of $\lambda.$ 
If $q=1$ then $\lambda$ is a star with the
black center so we will suppose that $q>1.$

Show first that if $\beta_2>1$ then there exists more than one tree of type
$\alpha$.
Indeed, let $v_1$ (resp. $v_2$) be a black vertex of $\lambda$ of valency
$\beta_1$ (resp. $\beta_2$) and let $f$ be a path connecting
$v_1$ and $v_2.$ Consider the following operation. Cut off a
branch $b$ of $\lambda$ growing from the vertex $v_2$
(that is a maximal subtree $b$ of $\lambda$ for which $v_2$
is a vertex of valency 1) such that $b$ does not contain
$f.$ Then glue $b$ to the vertex $v_1$ (see Fig. 9). 
\vskip 0.2cm
\medskip
\epsfxsize=9truecm
\centerline{\epsffile{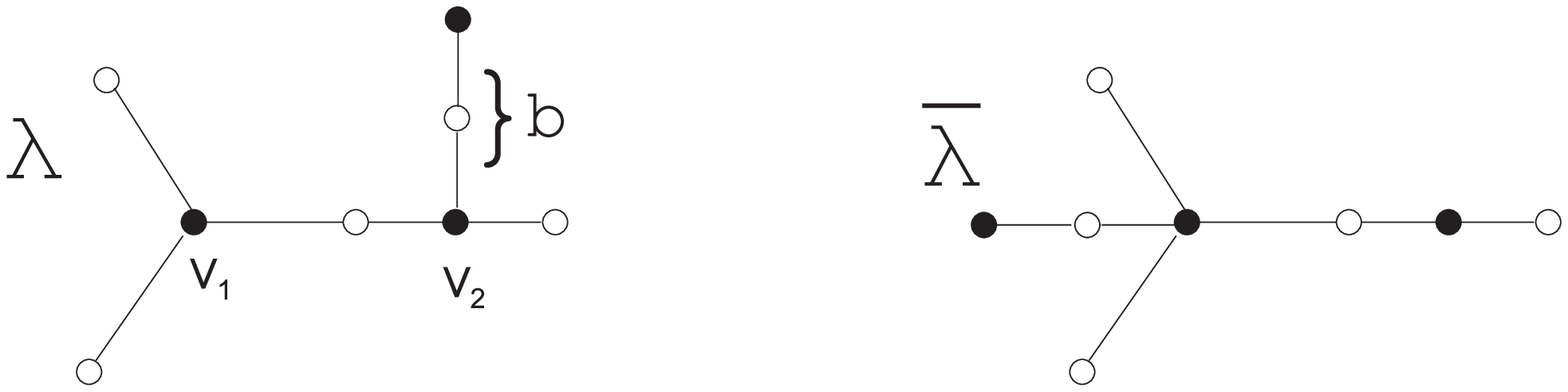}}
\smallskip
\centerline{Figure 9.}
\vskip 0.2cm
Since $\beta_2>1$ we always can perform such an operation 
and though generally there are many ways for doing it 
in any case the obtained tree
$\bar \lambda$
can not be isotopically equivalent to $\lambda$ since the maximal valency of
a black vertex of $\bar \lambda$ is greater than the corresponding
valency of $\lambda.$ On the other hand, the list of valencies of white vertices of $\bar \lambda$ remains the same. So, in the following
we will suppose that $\beta_2=1.$

If $\beta_1>2$ then cutting off any
branch $b$ of $\lambda$ growing from $v_1$
which does not contain
$f$ and gluing $b$ to the vertex $v_2$ we again obtain a tree
$\bar \lambda$ which is not isotopically equivalent to
$\lambda$ since the corresponding lists of valencies of black vertices $\{\beta_1,1,
1, ... ,1\}$ and $\{\beta_1-1,2,
1, ... ,1\}$ can not coincide due to the condition $\beta_1>2.$

Therefore, either $\beta_1=2$ and then $\lambda=\lambda_{r,s},$ or 
$\beta_1=1$ and then $\lambda$ is a star with
white center.
\vskip 0.2cm


\vskip 0.2cm

\noindent{\bf Acknowledgments.} I am grateful to G. Shabat
and A. Eremenko for discussions.
Also I would like to
thank
the Max Planck Institut f\"{u}r
Mathematik for its  support and hospitality.

\bibliographystyle{amsplain}

\begin{thebibliography}{10}

\bibitem {be} G. Belyi, \textit{
On Galois extensions of a maximal cyclotomic field,}
Math. USSR, Izv. 14, 1980, 247-256.


\bibitem {dh} A. Douady, J. Hubbard, \textit{
A proof of Thurston's topological characterization of rational functions}
Acta Math. 171, No.2, 1993, 263-297.

\bibitem {k} A. Kostrikin, \textit{Conservative polynomials}, in
``Stud. Algebra Tbilisi'', 115-129, 1984.

\bibitem {lz} S. Lando, A. Zvonkin, \textit{
Graphs on Surfaces and Their Applications,} 
Encyclopedia of Mathematical Sciences 141(II), Berlin: Springer, 2004.


\bibitem {pil} K. Pilgrim, \textit{Dessins d'enfants and Hubbard trees,}
Ann. Sci. \'Ecole Norm. Sup. (4) 33 (2000), no. 5, 671--693. 


\bibitem {poi1} A. Poirier, \textit{On postcritically finite polynomials,
part 1: critical portraits,} preprint, arxiv:math. DS/9305207.

\bibitem {poi2} A. Poirier, \textit{On postcritically finite polynomials,
part 2: Hubbard trees,} preprint, arxiv:math. DS/9307235.


\bibitem {sh} L. Schneps, \textit{Dessins d'enfants on the Riemann 
sphere}, in ``The Gro\-thendieck
Theory of Dessins D'enfants" (L. Shneps eds.), Cambridge University Press
``London mathematical society lecture notes series", vol. 200, 1994, 47-77.


\bibitem {sil} J. Silverman, \textit{
The field of definition for dynamical systems on $\P^1$,}
Compos. Math. 98, No.3, 269-304, 1995.

\bibitem {s} S. Smale, \textit{The fundamental theorem of algebra and
complexity theory,} Bull. Amer. Math. Soc., 4, 1-36, 1981.


\bibitem {t} D. Tischler, \textit{Critical points and values of complex 
polynomials,} J. of Complexity, 5, 438-456, 1989.




\end{thebibliography}

\end{document}